\documentclass[12pt]{article}
\usepackage{amsmath,amssymb,amscd,amsthm}

\makeatletter
\@namedef{subjclassname@2010}{%
  \textup{2010} Mathematics Subject Classification}
\makeatother

\theoremstyle{plain}
\newtheorem{theorem}{Theorem}[section]
\newtheorem{lemma}[theorem]{Lemma}
\newtheorem{coro}[theorem]{Corollary}
\newtheorem{propo}[theorem]{Proposition}

\theoremstyle{definition}

\theoremstyle{remark}

\renewcommand{\proof}{\noindent {\bf Proof.\ }}

\newcommand{\R}{\mathbb{R}}  

\newcommand{\mcal}[1]{\mathcal{#1}}

\newcommand{\supp}{\operatorname{supp}} 
\newcommand{\cardinal}[1]{\lvert #1\rvert} 

\newcommand{\st}{\hspace{0.1cm}:\hspace{0.1cm}} 
\newcommand{\dist}{\operatorname{dist}} 

\newcommand{\abs}[1]{\lvert #1 \rvert} 
\newcommand{\bigabs}[1]{\bigl\lvert #1 \bigr\rvert} 
\newcommand{\Bigabs}[1]{\Bigl\lvert #1 \Bigr\rvert} 

\newcommand{\linspan}{\operatorname{span}} 

\newcommand{\ip}[2]{\langle {#1},{#2}\rangle} 

\newcommand{\norm}[1]{\lVert #1 \rVert} 
\newcommand{\euclidnorm}[1]{\lvert #1 \rvert} 
\newcommand{\PP}{\operatorname{\mathbb{P}}} 
\newcommand{\EE}{\operatorname{\mathbb{E}}} 
\newcommand{\Var}{\operatorname{Var}} 
\newcommand{\E}{\mathbb{E}}

\begin{document}  

\title
{Smallest singular value of sparse random matrices}

\author{Alexander E. Litvak${}^{1}$ \qquad Omar Rivasplata}

\date{}

\maketitle

\footnotetext[1]{Research partially supported by  the
E.W.R. Steacie Memorial Fellowship.}

\begin{abstract}
We extend probability estimates on the smallest singular value of random matrices 
with independent entries to a class of sparse random matrices. We show that one can 
relax a previously used condition of uniform boundedness of the variances from below. 
This allows us to consider matrices with null entries or, more generally, with entries 
having small variances. Our results do not assume identical distribution of the entries 
of a random matrix and help to clarify the role of the variances of the entries. We also 
show that it is enough to require boundedness from above of the $r$-th moment, 
$r > 2$, of the corresponding entries.
\end{abstract}

\noindent
{\small \bf AMS 2010 Classification:}
{\small 46B06, 60B20, 15B52}


\noindent
{\small \bf Keywords:} {\small Random matrices, sparse matrices,
singular numbers, invertibility of random matrices, subgaussian
random variables, compressible and incompressible vectors,
deviation inequalities. }

\section{Introduction and main results}

Let $N \geq n$ be positive integers. 
In this paper we study the smallest singular value of $N \times n$ matrices 
$\Gamma = (\xi_{ji})$, whose entries are real-valued random variables 
obeying certain probability laws, and furthermore we are interested in allowing 
these matrices to contain some null entries (or, more generally, to contain 
entries with small variances).
Thus we deal with \emph{sparse} (or \emph{dilute}) random matrices.
Sparse random matrices and sparse structures play an important role,
as they arise naturally in many branches of pure and applied mathematics.
We refer to Chapter 7 of \cite{BSbook} for definitions,
relevant discussions, and references (see also the recent works \cite{GT, TV5}).

Understanding the properties of random matrices, in particular 
the behavior of their singular values (see the definitions in Section~\ref{prel}), 
is of importance in several fields, including Asymptotic Geometric Analysis, 
Approximation Theory, Probability and Statistics. 
The study of extreme singular values in classical random matrix
theory concentrates on their limiting behavior as the dimension grows to infinity. 
Such limiting behavior is now well understood for various kinds of random matrices 
whose entries are independent in aggregate, or independent up to the symmetry
constraints (e.g. hermitian or unitary matrices), in many cases even with
identical distribution being required.
We refer to the following books, surveys, and recent papers for history, results, 
and open problems in this direction
\cite{AGZ, BSbook, Erd, FS, Gir, mehta, TV5, TV4}.

In the \textit{non-limiting asymptotic} case very little was known till very recently.
In such a case one studies the rate of convergence, deviation inequalities, and the
general asymptotic behavior of singular values of a matrix as functions of the dimensions, 
assuming that the dimensions are large enough (growing to infinity).
The  Gaussian case, i.e. the case when the entries of the matrix are independent
$\mcal{N}(0,1)$ Gaussian, was treated independently in \cite{E} and \cite{Sz}
(see also \cite{Go} for related results, and the survey \cite{DS}).
In the last decade the attention shifted to other models, like matrices with independent
subgaussian entries (in particular, symmetric Bernoulli $\pm 1$ entries), independent
entries satisfying some moment conditions as well as matrices with independent columns
or rows satisfying some natural restrictions. Major achievements were obtained in
\cite{AGLPT, ALPT, lprt, R, RV1, RV2, tao-vu-2, tao-vu-3}.

In all previous non-limiting asymptotic results for random matrices with
independent entries, an important assumption was that the variances of all 
the entries are bounded below by one,  i.e. in a sense, that all entries are 
buffered away from zero and thus cannot be too small. 
Such a condition is not  natural for some 
applications, for instance when one deals with models in the theory of wireless
communications, where signals may be lost (or some small noise may appear), 
or with models in neural network theory, where the neurons are not of full
connectivity with each other, making sparse random matrices more suited in
modelling such partially connected systems. 

The main goal of our paper is to show that one
can significantly relax the condition of boundedness from below of all entries,
replacing it by averaging type conditions.
Thus our paper clarifies the role of the variances
in the corresponding previous results (cf. e.g. \cite{lprt, R, RV1, RV2}).
%
%
Another advantage of our results is that we require only boundedness
(from above) of the $r$-th moments for an arbitrary (fixed) $r > 2$.
We would like to emphasize that we don't require identical distributions
of all entries of a random matrix nor boundedness of the subgaussian moment 
of entries (both conditions were crucial for deep results of  \cite{RV2}). 
Moreover, the condition on entries ``to be identically distributed" is clearly 
inconsistent with our model, as, under such a condition, if one entry is zero
then automatically all entries are zeros.

\smallskip

We describe now our setting and results.
Our main results present estimates for the smallest singular value
$s_n(\Gamma)$ of large matrices $\Gamma$ of the type described.
It turns out the methods used to establish those estimates depend 
on the \emph{aspect ratio} of the matrices.
The aspect ratio of an $N \times n$ matrix $A$ is the ratio $n/N$
of number of columns to number of rows, or, more intuitively, the ratio
``width by height''.
To have a suggestive terminology, we will say that such matrix $A$ is
\begin{itemize}
\item ``\emph{tall}'' if $\frac{n}{N} \leq c_0$ for a small positive constant $c_0$;
\item ``\emph{almost square}'' if $\frac{n}{N}$ is close to 1.
\end{itemize}
Clearly, a matrix is square when its aspect ratio is equal to 1.

\smallskip

Since we will deal with random matrices under various conditions, 
for the sake of exposition clarity we list now all our conditions. 
For parameters $r > 2$, $\mu \geq 1$, $a_1>0$, $a_2>0$, $a_3 \in (0, \mu)$,
and $a_4 \in (0, 1]$, we will consider $N \times n$ random matrices
$\Gamma = (\xi_{ji})_{j\leq N,\, i\leq n}$ whose entries are \emph{independent}
real-valued \emph{centered} random variables satisfying the following conditions:
\begin{enumerate}
\item[(i)]   Moments: \quad $\EE |\xi_{ji}|^r \leq \mu^r$ \quad for all $j$ and $i$.
\item[(ii)]  Norm:  \quad $\PP \Bigl( \norm{\Gamma} > a_1 \sqrt{N} \Bigr) \leq e^{-a_{2} N}$.
\item[(iii)] Columns: \quad $\EE \norm{ (\xi_{ji})_{j=1}^N }_2^2 =
    \sum_{j=1}^{N} \EE \xi_{ji}^2 \geq a_3^2 N$ \quad for each $i$.
\end{enumerate}
For almost square and square matrices we also will need the following condition on rows.
\begin{enumerate}
\item[(iv)]  Rows:  \quad $\cardinal{ \{ i \st \EE \xi_{ji}^2 \geq 1 \} } \geq a_4 n$\quad for each $j$.
\end{enumerate}
Notice that these conditions allow our matrices to contain many null (or small) entries, 
in the sense that we don't impose any restrictions on the variance of a particular random variable. 
Naturally, in order for our random matrices to have entries of different kinds,
we do not require that the entries are identically distributed.
Our model is different from the sparse matrix models used e.g. in \cite{GT, TV5}, 
where zeros appeared randomly, i.e. starting from a random matrix whose entries have
variances bounded away from $0$, each entry was multiplied by another 
random variable of type $0/1$. 
Our model is more similar to those considered in \cite{Erd}, 
where a condition similar to (iii) was used for square symmetric matrices.

It is important to highlight that the parameters $\mu$, $r$, $a_1, a_2, a_3, a_4$
should be regarded as constants which do not depend on the dimensions $n$, $N$.
Note also that the ratio $\mu/a_3$ is of particular importance ($\mu$ is
responsible for the maximal $L_r$-norm of entries, while $a_3$ is an average-type  
substitution for the lower bound on $L_2$-norm of entries).

\smallskip


Before stating our main results let us comment our conditions in more detail.
The first condition is a standard requirement saying that the random variables are
not ``too big". For the limiting case it is known that one needs boundedness of
the forth moments. It turns out that for our estimates it is enough to ask
boundedness of moments of order $r= 2 + \varepsilon$ only, which improves
all previous results. In particular, this was one of the questions raised in
\cite{Vnew}, where the author proved corresponding estimates for entries with
bounded $4 + \varepsilon$ moment, and asked about $2 + \varepsilon$ moment.

The second condition is crucial for many results on random matrices. 
We recall that the norm of an $N \times n$ matrix is understood to be the operator norm 
from $\ell_2^n$ to $\ell_2^N$,
also called the spectral norm, which is equal to the largest singular value. In fact,
the question \emph{``What are the models of random matrices satisfying condition (ii)?"}
(and more generally, \emph{``What is the behavior of the largest singular value?"})
is one of the central questions in random matrix theory. Such estimates are well known 
for the Gaussian and subgaussian cases. We refer to \cite{ALPT, latala} and references therein 
for other models and recent developments on this problem.

We would like to emphasize that condition (ii) is needed in order to get probabilities
exponentially close to one. Alternatively, one may substitute this condition by
$$
  p_N := \PP \Bigl( \norm{\Gamma} > a_1 \sqrt{N} \Bigr) <1,
$$
in which case one should add $p_N$ to the estimates of probabilities in our theorems
below.

The main novelty in our model are conditions (iii) and (iv). 
These two conditions substitute the standard condition
\begin{equation}\label{old}
   \EE |\xi_{ji}|^2 \geq 1 \, \, \, \,  \mbox{ for all } \, \,  j,i,
\end{equation}
which was used in all previous works related to the smallest singular value  of
a random matrix (in the non-limiting case). 
Removing such strong assumption on \emph{all} entries, we allow
the possibility of zeros to appear among the entries of a random matrix.
Our conditions (iii) and (iv) 
should be compared with the normalization conditions (1.1) and (1.16) in \cite{Erd}.
Our methods are similar to those used in \cite{lprt, RV2}, 
however we deal with a rather different model,
and correspondingly our proofs require much more delicate computations. 
In particular, the proof of key Proposition~\ref{propo-tall}, which estimates the probability 
that for a fixed vector $x$ the Euclidean  norm $\|\Gamma x\|_2$ is small, 
is much more involved
(cf. the proof of \cite[Proposition 3.4]{lprt} or \cite[Corollary 2.7]{RV1}).

Of course we want to rule out matrices having a column or a row consisting of zeros only,
for if there is a zero column then immediately $s_n(\Gamma)=0$, while if there is a zero row
then the matrix $\Gamma$ is essentially of size $(N-1)\times n$.
Hence we need some general assumptions on the columns and the rows of the matrices under
consideration. Our condition (iii) alone  implies that each column vector of
the matrix has relatively big $\ell_2$-norm. Moreover, condition (iii)
together with condition (i)
guarantee that proportionally many rows have $\ell_2$-norms bounded away from 0.
It turns out that condition (iii) is already enough for ``tall'' matrices, when $N > Cn$,
as the first theorem below shows. The cases of ``almost square'' and square matrices are
more delicate, because $N$ becomes closer to $n$, and we need to control the behavior of
rows more carefully. Condition (iv) ensures that each row of the matrix
has proportionally many entries with variance at least one.

\medskip


Now we state our results. The first theorem deals with ``tall'' matrices
and extends the corresponding result from \cite{lprt} 
(for uniformly \textit{bounded above} mean zero random variables 
with bounded below variances this was shown in \cite{BDG}). 
Note that we use only three conditions, (i), (ii), and (iii), 
while condition (iv) is not required for this result.

\begin{theorem}
\label{teo-tall}
Let $r > 2$, $\mu \geq 1$, $a_1, a_2, a_3 > 0$ with $ a_3 < \mu$.
Let $1 \leq n < N$ be integers, and write $N$ in the form $N = (1 + \delta)n$.
Suppose $\Gamma$ is an $N \times n$ matrix whose entries are independent
centered random variables
such that conditions (i), (ii) and (iii) are satisfied.
There exist positive constants $c_1$, $c_2$ and $\delta_0$
(depending only on the parameters $r$, $\mu$, $a_1$, $a_2$, $a_3$)
such that whenever $\delta \geq \delta_0$, then
$$
\PP \Bigl( s_n(\Gamma) \leq c_1 \sqrt{N} \Bigr) \leq e^{-c_2 N}.
$$
\end{theorem}

\noindent
{\bf Remark. } Our proof gives that $c_1 = c_1 (r, \mu, a_3)$,
$c_2 = c_2 (r, \mu, a_2, a_3)$ and $\delta_0 = \delta_0 (r, \mu, a_1, a_3)$.

\medskip

Our next theorem is about ``almost square'' matrices.
This theorem extends
\cite[Theorem 3.1]{lprt}. Here both conditions (iii) and (iv) are needed
in order to substitute condition (\ref{old}).

\begin{theorem}
\label{teo-almost}
Let $r > 2$, $\mu \geq 1$, $a_1, a_2>0$, $a_3\in (0, \mu)$, $a_4 \in (0, 1]$.
Let $1 \leq n < N$ be integers, and write $N$ in the form $N = (1 + \delta)n$.
Suppose $\Gamma$ is an $N \times n$ matrix whose entries are independent
centered random variables such that conditions (i), (ii), (iii) and (iv)
are satisfied. There exist positive constants  $c_1$, $c_2$,
$\tilde{c}_1$ and $\tilde{c}_2$,
depending only on the parameters $r$, $\mu$, $a_1$, $a_2$, $a_3$, $a_4$,
and a positive constant $\gamma =\gamma(r,\mu, a_1, a_3)< 1$,
such that if
$$
  a_4> 1 - \gamma  \hspace{1cm}
  \text{ and } \hspace{1cm} \delta \geq \frac{\tilde{c}_1}{ \ln(2+ \tilde{c}_2 n) }
$$
then
$$
  \PP \Bigl( s_n(\Gamma) \leq c_1 \sqrt{N} \Bigr) \leq e^{-c_2 N}.
$$
\end{theorem}

\noindent
{\bf Remarks. \ }
{\bf 1. }
Our proof gives that
$c_1 = c_1 (r,\mu,a_1,a_3,\delta)$,
$c_2 = c_2 (r,\mu,a_2,a_3)$, $\tilde{c}_1 = \tilde{c}_1 (r,\mu,a_1,a_3)$ and
$\tilde{c}_2 = \tilde{c}_2 (r, \mu, a_1, a_3, a_4)$. \\
{\bf 2. }
Note that for small $n$, say for $n \leq 2/\tilde{c}_2$, Theorem~\ref{teo-almost}
is trivial for every $\delta>0$, either by adjusting the constant $c_2$ (for small $N$)
or by using Theorem~\ref{teo-tall} (for large $N$).

\smallskip

Let us note that in a sense our Theorems~\ref{teo-tall} and \ref{teo-almost}
are incomparable with the corresponding result of \cite{RV2}. 
First, we don't restrict our results only to the subgaussian case.
The requirement of boundedness of the subgaussian moment is much
stronger, implying in particular boundedness of moments of all orders,
which naturally yields stronger estimates.
Second, another condition essentially used in \cite{RV2} is
``entries are identically distributed."
As was mentioned above, such a condition is inconsistent with our model,
since having one zero we immediately get the zero matrix.

\medskip

Our third theorem shows that we can also extend to our setting
the corresponding results from \cite{RV1}, where the i.i.d. case was treated, and 
from \cite{AGLPT1, AGLPT}, which dealt with the case of independent log-concave columns.
Note again that we work under the assumption of bounded $r$-th moment 
(for a fixed $r > 2$).
In fact in \cite{RV1} two theorems about square matrices were proved. 
The first one is for random matrices whose entries have bounded fourth moment. 
Our Theorem~\ref{teo-square} extends this result with much better probability. 
The second main result of \cite{RV1} requires the boundedness of subgaussian 
moments as well as identical distributions of entries in each column, and, thus, 
is incomparable with Theorem~\ref{teo-square}.

\begin{theorem}
\label{teo-square}
Let $r > 2$, $\mu \geq 1$, $a_1, a_2, a_3, a_4 > 0$ with $ a_3 < \mu$.
Suppose $\Gamma$ is
an $n \times n$ matrix whose entries are independent centered random variables
such that conditions (i), (ii), (iii) and (iv) are satisfied. 
Then there exists a positive constant $\gamma_0 =\gamma_0(r,\mu, a_1, a_3)< 1$ 
such that if $a_4> 1 - \gamma_0$ then for every $\varepsilon \geq 0$
\begin{displaymath}
 \PP \bigl( s_n (\Gamma) \leq \varepsilon n^{-1/2} \bigr) \leq C \bigl(
 \varepsilon + n^{1-r/2} \bigr),
\end{displaymath}
where $C$ depends on the parameters $r$, $\mu, a_1, a_2, a_3, a_4$.
\end{theorem}

\smallskip

Finally we would like to mention that all results can be extended to the
complex case in a standard way.

\bigskip

\noindent
{\bf Acknowledgment.}
The authors would like to thank N. Tomczak-Jaegermann for many useful conversations.
We also thank S.~Spektor for showing us reference \cite{petrov}
and S. O'Rourke for showing us reference \cite{Erd}.
The second named author thanks G.~Schechtman for hosting him at the
Weizmann Institute of Science in Spring 2008, during which time part of this work
was done.

\section{Notation and preliminaries}\label{prel}

We start this section by agreeing on the notation that we will use throughout.
For $1 \leq p \leq \infty$, we write $\norm{x}_p$ for the $\ell_p$-norm of
$x = (x_i)_{i \geq 1}$, i.e. the norm defined by
\begin{displaymath}
\norm{x}_p = \Bigl( \sum_{i\geq 1} \abs{x_i}^p \Bigr)^{1/p} \text{ for } p < \infty \qquad
\text{ and } \qquad \norm{x}_{\infty} = \sup_{i \geq 1} \abs{x_i}.
\end{displaymath}
Then, as usual, $\ell_p^{n} = ( \R^n, \norm{\cdot}_p )$.
The unit ball of $\ell_p^{n}$ is denoted $B_p^n$.
Also, $S^{n-1}$ denotes the unit sphere of $\ell_2^n$, and $e_1,\ldots,e_n$ is the canonical
basis of $\ell_2^n$.

We write $\ip{\cdot}{\cdot}$ for the standard inner product on $\R^n$.
By $\euclidnorm{x}$ we denote the standard Euclidean norm (i.e. $\ell_2$-norm)
of the vector $x = (x_i)_{i \geq 1}$.
On the other hand, when $A$ is a set, by $\cardinal{A}$ we denote the cardinality of $A$.

The support of a vector $x = (x_i)_{i \geq 1}$, meaning the set of
indices corresponding to nonzero coordinates of $x$, is denoted by $\supp(x)$.

Given a subspace $E$ of $\R^n$ we denote by $P_{E}$ the orthogonal projection onto $E$.
If $E = \R^{\sigma}$ is the coordinate subspace corresponding to a set of
coordinates $\sigma \subset \{ 1,\ldots,n \}$, we will write $P_{\sigma}$ as a shorthand
for $P_{\R^{\sigma}}$.

\medskip

Let $\mcal{N}\subset D\subset \R^n$ and  $\varepsilon > 0$. 
Recall that $\mcal{N}$ is called an $\varepsilon$-net of $D$ (in the Euclidean metric) if
$$
  D \subset \bigcup_{v \in \mcal{N}} (v + \varepsilon B_2^n).
$$
In case $D$ is the unit sphere $S^{n-1}$ or the unit ball $B_2^{n}$, a well known
volumetric argument (see for instance \cite[Lemma 2.6]{ms}) establishes that
for each $\varepsilon > 0$ there is an $\varepsilon$-net $\mcal{N}$ of $D$ with cardinality
$\cardinal{\mcal{N}} \leq (1 + 2/\varepsilon)^n$.

\medskip

\subsection{Singular values.}

Suppose $\Gamma$ is an $N \times n$ matrix with real entries.
The singular values of $\Gamma$, denoted $s_k(\Gamma)$, are the eigenvalues of the
$n \times n$ matrix $\sqrt{\Gamma^t \,\Gamma}$, arranged in the decreasing order.
It is immediate that the singular values are all non-negative, and further the number
of nonzero singular values of $\Gamma$ equals the rank of $\Gamma$.

The largest singular value $s_1(\Gamma)$ 
and the smallest singular value $s_n(\Gamma)$
are particularly important.
They may be equivalently given by the expressions
\begin{equation*}
s_1 (\Gamma) 
= \|\Gamma: \ \ell_2^n \to \ell_2^N\|
= \sup \bigl\{ \euclidnorm{\Gamma x} \st \euclidnorm{x} = 1 \bigr\},
\quad 
s_n (\Gamma) 
= \inf \bigl\{ \euclidnorm{\Gamma x} \st \euclidnorm{x} = 1 \bigr\}.
\end{equation*}
In particular for every vector $x\in \mathbb{R}^n$ one  has
\begin{equation}\label{norm-estimate}
s_{n} (\Gamma) \euclidnorm{x} \leq \euclidnorm{\Gamma x} \leq s_1 (\Gamma) \euclidnorm{x}.
\end{equation}

Note that the estimate on the left-hand side becomes trivial if $s_n (\Gamma) = 0$.
On the other hand, when $s_n (\Gamma) > 0$ the matrix $\Gamma$ is a bijection 
on its image, and can be regarded as an embedding from $\ell_2^n$ into $\ell_2^N$, 
with \eqref{norm-estimate} providing an estimate for the distortion of the norms 
under $\Gamma$.

To estimate the smallest singular number, we will be using the following equivalence,
which clearly holds for every matrix $\Gamma$ and every  $\lambda \geq 0$:
\begin{equation} \label{smallest-sv-equivalence}
s_n (\Gamma) \leq
\lambda \quad \Longleftrightarrow \quad
\exists x \in S^{n-1} \st \euclidnorm{\Gamma x} \leq \lambda.
\end{equation}

\medskip

\subsection{Subgaussian random variables.}

All random quantities appearing in this work are defined on 
the same underlying probability space $(\Omega,\mcal{A},\PP)$.
We will present estimates for the smallest singular value of matrices whose entries
are independent random variables satisfying certain assumptions.
Our results are valid for a large class of matrices which includes, in particular,
those whose entries are subgaussian random variables.

A (real-valued) random variable $X$ is called \emph{subgaussian} when there exists 
a positive constant $b$ such that for every $t \in \R$
\begin{equation*}
\EE e^{t X} \leq e^{b^2 t^2 / 2}.
\end{equation*}
When this condition is satisfied with a particular value of $b > 0$, we also say that
$X$ is $b$-subgaussian, or subgaussian with parameter $b$.
The minimal $b$ in this capacity is called the \emph{subgaussian moment} of $X$.

It is an easy consequence of this definition that if $X$ is $b$-subgaussian,
then $\EE(X) = 0$ and $\Var(X) \leq b^2$.
Thus all subgaussian random variables are centered.
The next proposition
presents well-known equivalent conditions for a centered random
variable to be subgaussian.

\begin{propo}
For a centered random variable $X$, the following statements are equivalent:
\begin{itemize}
\item[(1)] $\exists b > 0, \hspace{2mm} \forall t \in \R, \hspace{2mm}
\EE e^{t X} \leq e^{b^2 t^2 / 2} $
\item[(2)] $\exists b > 0, \hspace{2mm} \forall \lambda > 0, \hspace{2mm}
\PP (\abs{X} \geq \lambda) \leq 2 e^{-\lambda^2 /  b^2} $
\item[(3)] $\exists b > 0, \hspace{2mm} \forall p \geq 1, \hspace{2mm}
(\EE \abs{X}^p)^{1/p} \leq  b  \sqrt{\smash[b]p}$
\item[(4)] $\exists c > 0, \hspace{2mm} \EE e^{c X^2} < +\infty $
\end{itemize}
\end{propo}

Two important examples of subgaussian random variables are the centered Gaussian
themselves and the symmetric Bernoulli  $\pm 1$ random variables.
In general, any centered and bounded random variable is subgaussian.

We point out that, as consequence of the subgaussian tail estimate, the norm of a
matrix whose entries are independent subgaussian random variables is of the order
of $\sqrt{N}$ with high probability.
Namely, the following proposition holds (see e.g. \cite[Fact 2.4]{lprt},  where
this was shown for symmetric random variables, the case of centered is essentially
the same).

\begin{propo}
\label{norm-subgaussian}
Let $N \geq n \geq 1$ be positive integers.
Suppose $\Gamma$ is an $N \times n$ matrix whose entries are independent subgaussian 
random variables with subgaussian parameters bounded above uniformly by $b$.
Then there are positive constants $c, C$ (depending only on $b$) such that for every $t>C$
$$
  \PP \bigl( \norm{\Gamma} > t \sqrt{N} \bigr) \leq e^{-c t^2 N}.
$$
\end{propo}

\subsection{Compressible and incompressible vectors.}

As equivalence \eqref{smallest-sv-equivalence} suggests,
to estimate the smallest singular value of $\Gamma$ we estimate
the norm $\euclidnorm{\Gamma x}$ for vectors $x \in S^{n-1}$.
More precisely, we will estimate $\euclidnorm{\Gamma x}$ individually
for vectors in an appropriately chosen $\varepsilon$-net and,
as usual, we use the union bound.
In the case of ``tall'' matrices just one single $\varepsilon$-net is enough
for this approximation method to work; but in the case of ``almost square'' matrices,
as well as for square matrices, we will need to split the sphere into two parts
according to whether the vector $x$ is compressible or incompressible,
in the sense that we now define.

\medskip

Let $m \leq n$ and $\rho \in (0,1)$. A vector $x \in \R^n$ is called
\begin{itemize}
\item \emph{$m$-sparse} 
          if $\cardinal{\supp(x)} \leq m$, that is, if $x$ has at most $m$ nonzero entries.
\item \emph{$(m,\rho)$-compressible} 
          if it is within Euclidean distance $\rho$ from the set of all $m$-sparse vectors.
\item \emph{$(m,\rho)$-incompressible} 
          if it is not $(m,\rho)$-compressible.
\end{itemize}
The sets of sparse, compressible, and incompressible vectors will be denoted, respectively,
$Sparse(m)$, $Comp(m,\rho)$, and $Incomp(m,\rho)$.
The idea to split the Euclidean sphere into two parts
goes back to Kashin's work \cite{Kash} on orthogonal decomposition of $\ell_1^{2n}$,
where the splitting was defined using the ratio of $\ell _2$ and $\ell_1$ norms.
This idea was recently used by Schechtman (\cite{Sch}) in the same context.
The splitting the sphere essentially as described above
appeared in \cite{lprt, lprtvcr} and was later used in many works
(e.g. in \cite{RV1, RV2}).

\medskip

It is clear from these definitions that, for a vector $x$, the following holds:
$$
   x \in Comp(m,\rho) \ \   \Longleftrightarrow \ \
   \exists \sigma \subset \{ 1,\ldots,n \} \ \text{ with }\ \cardinal{\sigma^c}
   \leq m \ \text{ such that }\  \euclidnorm{P_{\sigma} x} \leq \rho \notag\\
$$
\begin{equation}
\label{incompressible-equivalence}
x \in Incomp(m,\rho) \ \  \Longleftrightarrow \ \
\forall \sigma \subset \{ 1,\ldots,n \} \  \text{ with }\ \cardinal{\sigma^c} \leq m
\ \text{ one has }\ \euclidnorm{P_{\sigma} x} > \rho.
\end{equation}

\subsection{Two more results.}

Here we formulate two results, which will be used in the next section.
The first one is a quantitative version of the Central Limit Theorem (CLT), 
called Berry-Ess\'{e}en inequality.
The second one is a general form of the Paley-Zygmund inequality 
(see e.g. \cite[Lemma 3.5]{lprt}).

\begin{theorem}[Berry-Ess\'{e}en CLT]
\label{berry-esseen} Let $2 < r \leq 3$.
Let $\zeta_1, \ldots, \zeta_n$ be independent centered random variables
with finite $r$-th moments and set
$\sigma^2 := \sum_{k = 1}^{n} \EE \abs{\zeta_k}^2$.
Then for all $t \in \R$
\begin{displaymath}
 \biggl| \PP \biggl( \frac{1}{\sigma}\sum_{k = 1}^{n} \zeta_k \leq t \biggr) -
 \PP \bigl( g \leq t \bigr) \biggr| \leq \frac{C}{\sigma^r} \sum_{k = 1}^{n}
 \EE \abs{\zeta_k}^r,
\end{displaymath}
where $g \sim \mcal{N}(0,1)$ and $C$ is an absolute constant.
\end{theorem}

\medskip

\noindent{\bf Remarks.}
\newline
{\bf 1.\ } The standard form of Berry-Ess\'{e}en inequality requires 
finite 3-rd moment (i.e., it is usually stated for $r = 3$),
see e.g. \cite[p. 544]{fellerII} or \cite[p. 300]{loeveI}.
The form used here is from \cite{petrov} (see Theorem~5.7 there).
\newline
{\bf 2.\ } If $r \geq 3$, then clearly we have boundedness 
of $3$-rd moment for free, and in this case we use 
the standard form of Berry-Ess\'{e}en inequality (i.e., with $r = 3$).

\medskip

\begin{lemma}[Paley-Zygmund inequality]
\label{pa-zy}
Let $p \in (1, \infty)$,  $q = p/(p-1)$.
Let $f \ge 0$ be a random variable with $\EE f^{2p} < \infty$. Then
for every  $0\leq \lambda \leq \sqrt{\EE f^2}$ we have
$$
 \PP \left( f > \lambda  \right)  \ge \frac{(\EE f^2 - \lambda^2)^q}{(\EE f^{2p})^{q/p}}.
$$
\end{lemma}

\section{Small ball probabilities for random sums}

In this section we gather auxiliary results related to random sums, their small ball
probabilities, etc., which are needed later. In fact, we adjust corresponding results
from \cite{lprt} and \cite{RV1} to our setting. These results are also of independent interest. 
We provide proofs for the sake of completeness.

The following lemma provides a lower bound on the small ball probability of a random sum.
Its proof follows the steps of \cite[Lemma~3.6]{lprt} with the appropriate
modification to deal with centered random variables (rather than symmetric),
to remove the assumption that the variances are bounded from below uniformly,
and to replace the condition of finite 3-rd moments by finite $r$-th moments ($r > 2$).

\begin{lemma}
\label{lema-deviation}
Let $2<r \leq 3$ and $\mu \geq 1$.
Suppose $\xi_1, \ldots, \xi_n$ are independent centered random variables
such that 
$\EE \abs{ \xi_{i} }^r \leq \mu^r$ for every $i = 1, \ldots, n$.
Let $x = (x_i) \in \ell_{2}$ be such that $\euclidnorm{x} = 1$.
Then for every $\lambda \geq 0$
\begin{equation*}
\PP \biggl( \Bigabs{ \sum_{i=1}^{n} \xi_i x_i } > \lambda \biggr) \geq
\left( \frac{ [\EE \sum_{i = 1}^{n} \xi_i^2 x_i^2 - \lambda^2]_+ }{ 8\mu^2 } \right)^{r/(r-2)}.
\end{equation*}
\end{lemma}

\proof
Define $f = \bigabs{\sum_{i=1}^{n} \xi_i x_i}$.
Let $\varepsilon_1,\ldots,\varepsilon_n$ be independent symmetric Bernoulli
$\pm 1$ random variables,
which are also independent of $\xi_1\dots,\xi_n$.
Using the symmetrization inequality \cite[Lemma 6.3]{ledoux-talagrand},
and applying Khinchine's inequality, 
we obtain
\begin{align*}
\EE f^r
\leq 2^r \EE \Bigabs{ \sum_{i=1}^{n} \varepsilon_i \xi_i x_i }^r
= 2^r \EE_{\xi} \EE_{\varepsilon} \Bigabs{ \sum_{i \geq 1} \varepsilon_i \xi_i x_i }^r
\leq 2^r 2^{r/2} \EE_{\xi} \biggl( \sum_{i \geq 1} \xi_i^2 x_i^2 \biggr)^{r/2}.
\end{align*}

Now consider the set
$$
   \mcal{S} := \biggl\{ s=(s_i) \in \ell_1 \st s_i \geq 0\;
  \text{ for every $i$ and } \sum_{i \geq 1} s_i = 1 \biggr\}.
$$
We define a function $\varphi : \mcal{S} \to \R$ by
$$
\varphi(s) = \EE_{\xi} \biggl( \sum_{i \geq 1} \xi_i^2 s_i \biggr)^{r/2}.
$$
This function is clearly convex, so that
$$
\sup_{s \in \mcal{S}} \varphi(s)
= \sup_{i\geq 1} \varphi(e_i)
= \sup_{i\geq 1} \EE_{\xi} (\xi_i^2)^{r/2}
\leq \mu^r.
$$
Thus $\EE f^r \leq 2^{3r/2} \mu^r$.
On the other hand, using the independence of $\xi_1, \ldots, \xi_n$,
$$
  \EE f^2 = \EE \sum_{i \geq 1} \xi_i^2 x_i^2.
$$
Lemma~\ref{pa-zy} with $p = r/2$, $q = r/(r-2)$ implies the desired estimate.
\qed

\medskip

The next proposition, which is a  consequence of
Theorem~\ref{berry-esseen},  allows us to estimate the small ball probability.
The proof goes along the same lines as the proof of \cite[Proposition 3.2]{lprt}
(see also \cite[Proposition 3.4]{lprtv}),
with slight modifications to remove the assumption about variances.
Recall that for a subset $\sigma \subset \{ 1, 2, \ldots, n \}$,
 $P_{\sigma}$ denotes the coordinate projection onto $\R^{\sigma}$.

\begin{propo}
\label{propo-sbp}
Let $2<r \leq 3$ and $\mu \geq 1$.
Let $(\xi_i)_{i=1}^n$ be independent centered random variables with
$\E |\xi_{i}|^r \leq \mu^r$
for all $i = 1, 2, \ldots, n$.
There is a universal constant $c > 0$ such that
\begin{enumerate}
\item[(a)] For every  $a < b$ and every $x = (x_i) \in \R^n$ satisfying
\mbox{$A :=  \sqrt{ \EE \sum_{i=1}^{n} \xi_i^2 x_i^2} >0$}
one has
$$
 \PP \biggl( a \leq \sum_{i=1}^{n} \xi_i x_i < b \biggr)
 \leq \frac{b-a}{\sqrt{2 \pi} A} + c \biggl( \frac{ \norm{x}_r }{A}
 \mu \biggr)^r .
$$
\item[(b)] For every $t > 0$, every $x = (x_i) \in \R^n$ and every
$\sigma \subset \{ 1, 2, \ldots, n \}$ satisfying
$A_{\sigma} := \sqrt{ \EE \sum_{i \in \sigma} \xi_i^2 x_i^2} >0$
one has
$$
  \sup_{v \in \R} \PP \biggl( \Bigabs{ \sum_{i=1}^{n} x_i \xi_i - v } <
  t \biggr)
   \leq \frac{2t}{\sqrt{2 \pi} A_{\sigma}} + c \biggl( \frac{
   \norm{P_{\sigma}x}_r }{A_{\sigma}} \mu \biggr)^r .
$$
\end{enumerate}
\end{propo}

The next corollary gives an estimate on the small ball probability
in the spirit of \cite[Corollary 2.10]{RV1}.

\begin{coro}
\label{coro-sbp-big-eps}
Let $2<r \leq 3$ and $\mu \geq 1$.
Let $\xi_1, \ldots, \xi_n$ be independent centered random variables
with $\E |\xi_{i}|^r \leq \mu^r$
for every $i = 1, \ldots, n$. Suppose $x = (x_i)\in \R^n$ and
$\sigma \subset \{ 1, \ldots, n \}$ are such that $A \leq \abs{x_i} \leq B$
and $\EE \xi_i^2 \geq 1$ for all $i \in \sigma$. Then  for all $t \geq 0$
\begin{displaymath}
  \sup_{v \in \R} \PP \biggl( \Bigabs{ \sum_{i=1}^{n} x_i \xi_i - v } <
  t \biggr) \leq \frac{C}{\cardinal{\sigma}^{r/2 - 1}} \biggl( \frac{t}{A} +
  \mu^r \Bigl( \frac{B}{A} \Bigr)^r \biggr),
\end{displaymath}
where $C > 0$ is an absolute constant.
\end{coro}

\proof
By assumptions on coordinates of $x$ we have
$$
  A_\sigma^2 := \EE \sum_{i \in \sigma} \xi_i^2 x_i^2 \geq \cardinal{\sigma} A^2
$$
and
$$
 \norm{P_{\sigma} x}_r^r = \sum_{i \in \sigma} \abs{x_i}^r
 \leq \cardinal{\sigma} B^r.
$$
Then, by part (b) of Proposition \ref{propo-sbp}
\begin{align*}
 \sup_{v \in \R} \PP \biggl( \Bigabs{ \sum_{i=1}^{n} x_i \xi_i - v } <
  t \biggr) &\leq
  \sqrt{ \frac{2}{\pi} } \frac{t}{A\cardinal{\sigma}^{1/2}} +
  c \mu^r \frac{ B^r \cardinal{\sigma} }{A^r \cardinal{\sigma}^{r/2}} \\
  &\leq \frac{C}{\cardinal{\sigma}^{r/2 - 1}} \biggl( \frac{t}{A} + \mu^r
  \Bigl( \frac{ B}{A} \Bigr)^r \biggr).
\end{align*}
\qed
\medskip

We need the following lemma proved in \cite[Lemma 3.4]{RV1}.

\begin{lemma} 
\label{incompressible-are-spread}
Let $\gamma,\rho \in (0,1)$, and let $x \in Incomp(\gamma n,\rho)$.
Then there exists a set $\sigma=\sigma_x \subset \{ 1, \ldots, n \}$ of
cardinality $\cardinal{\sigma} \geq \frac{1}{2} \rho^2 \gamma n$
and such that for all $k \in \sigma$
\begin{displaymath}
\frac{\rho}{\sqrt{2n}} \leq \abs{x_k} \leq \frac{1}{\sqrt{\gamma n}}.
\end{displaymath}
\end{lemma}

The next lemma is a version of \cite[Lemma 3.7]{RV1}, modified in order to
remove the assumption ``$variances \geq 1$''.

\begin{lemma}
\label{sbp-incompressible}
Let $2<r \leq 3$ and $\mu \geq 1$.
  Let $\xi_1, \ldots, \xi_n$ be independent centered random variables with
$\E |\xi_{i}|^r \leq \mu^r$
for every $i$. Suppose $\overline{\sigma} :=
  \{ i \st \EE \xi_i^2 \geq 1 \}$ has cardinality $\cardinal{ \overline{\sigma} } \geq a_4 n$.
  Let $\gamma, \rho \in (0,1)$, and consider a vector $x \in Incomp(\gamma n,\rho)$.
 Assuming that $a_4 + \frac{1}{2} \rho^2 \gamma >1$ we have for every $t \geq 0$
\begin{displaymath}
  \sup_{v \in \R} \PP \biggl( \Bigabs{ \sum_{i=1}^{n} x_i \xi_i - v }
  < t \biggr) \leq c ( t n^{\frac{3-r}{2}} + \mu^r n^{\frac{2-r}{2}} ),
\end{displaymath}
where $c$ is a positive constant which depends on $\gamma$, $\rho$, $a_4$, and $r$.
\end{lemma}

\proof
Let $\sigma_{x}$ be the set of spread coefficients of $x$ from
Lemma \ref{incompressible-are-spread}, so that
$\cardinal{\sigma_{x}} \geq \frac{1}{2}\rho^2 \gamma n$.
Set $\sigma := \overline{\sigma} \cap \sigma_{x}$. Then
\begin{displaymath}
\cardinal{\sigma} = \cardinal{\overline{\sigma}} + \cardinal{\sigma_{x}} -
\cardinal{ \overline{\sigma} \cup \sigma_{x} }
\geq a_4 n + \frac{1}{2}\rho^2 \gamma n - n =: c_0 n.
\end{displaymath}
By the construction, for every $i \in \sigma$ we have
$$
  \frac{\rho}{\sqrt{2n}} \leq \abs{x_i} \leq \frac{1}{\sqrt{\gamma n}}.
$$
Applying  Corollary \ref{coro-sbp-big-eps} we obtain
\begin{align*}
\sup_{v \in \R} \PP \biggl( \Bigabs{ \sum_{i=1}^{n} x_i \xi_i - v }
  < t \biggr)
& \leq \frac{C}{\cardinal{\sigma}^{r/2 - 1}}
\biggl( \frac{\sqrt{2n}t}{\rho} + \mu^r
\Bigl( \frac{\sqrt{2}}{\rho\sqrt{\gamma}} \Bigr)^r \biggr)\\
&\leq \frac{C}{(c_0 n)^{r/2 - 1}}
\biggl( \frac{\sqrt{2n}t}{\rho} + \mu^r
\Bigl( \frac{\sqrt{2}}{\rho\sqrt{\gamma}} \Bigr)^r \biggr)\\
&\leq c ( t n^{\frac{3-r}{2}} + \mu^r n^{\frac{2-r}{2}} ).
\end{align*}
\qed

\section{``Tall'' matrices (proof of Theorem~\ref{teo-tall})}

In this section we prove Theorem~\ref{teo-tall}, which establishes an estimate on the
smallest singular value for ``tall'' random matrices, meaning matrices whose aspect
ratio $n/N$ is bounded above by a small positive constant (independent of $n$ and $N$).
It is important to notice that Theorem~\ref{teo-tall} uses only conditions (i), (ii),
and (iii), i.e. no condition on the rows is required here.

The proof depends upon an estimate on the norm $\abs{ \Gamma x }$
for a fixed vector $x$, which is provided by the following proposition.

\begin{propo}
\label{propo-tall}
Let $1 \leq n < N$ be positive integers.
Suppose $\Gamma$ is a matrix of size $N \times n$ whose entries are independent centered
random variables satisfying conditions (i), (ii) and (iii) for some $2<r \leq 3$, $\mu \geq 1$
and $a_1, a_2, a_3 > 0$ with $a_3 < \mu$. Then for every $x \in S^{n-1}$ we have
$$
\PP \Bigl( \euclidnorm{\Gamma x} \leq b_{1} \sqrt{N} \Bigr)
\leq e^{-b_{2} N },
$$
where $b_{1},b_{2} > 0$ depend only on $\mu$, $a_{3}$ and $r$.
\end{propo}

\noindent
{\bf Remark. } Our proof gives that
$$
b_{1} = \frac{a_3^{4}}{2^{5}\mu^2} \Bigl( \frac{a_3^{2}}{2^{5}\mu^2} \Bigr)^{r/(r - 2)}, \qquad
b_{2} = \frac{a_3^{2}}{2^{3}\mu^2} \Bigl( \frac{a_3^{2}}{2^{5}\mu^2} \Bigr)^{r/(r - 2)}.
$$

We postpone the proof of this technical result to the last section,
so that we may keep the flow of our exposition uninterrupted.

\medskip

\noindent
{\bf Proof of Theorem \ref{teo-tall}.} \quad Passing to $r_0=\min\{3, r\}$ we may
assume without loss of generality that  $r\leq 3$.

Let $t\geq 0$ and $\Omega_{0} := \{ \omega \st \norm{\Gamma} \leq a_1 \sqrt{N} \}$.
By \eqref{smallest-sv-equivalence} it is enough to estimate
the probability of the event
$$
  E := \{ \omega \st \exists x \in S^{n-1} 
\hspace{1mm} \text{ s.t. } \hspace{1mm} \euclidnorm{\Gamma x} \leq t \sqrt{N} \}.
$$
To this end we use the inclusion $E \subset (E \cap \Omega_{0}) \cup \Omega_{0}^c$
and the union bound.

To estimate $\PP (E\cap \Omega_{0})$,
let $0 < \varepsilon \leq 1$, and let $\mcal{N}$ be an $\varepsilon$-net of $S^{n-1}$ 
with cardinality $\cardinal{ \mcal{N} } \leq (3/\varepsilon)^n$.
For any $x \in S^{n-1}$ we can find $y \in \mcal{N}$ such that $\euclidnorm{x - y} \leq \varepsilon$.
If further $x$ satisfies $\euclidnorm{\Gamma x} \leq t \sqrt{N}$, then the corresponding $y$
satisfies
\begin{equation}\label{normy}
 \euclidnorm{\Gamma y} \leq \euclidnorm{\Gamma x} + \norm{\Gamma}\cdot\euclidnorm{y - x}
 \leq t \sqrt{N} + \varepsilon a_1 \sqrt{N} = (t + \varepsilon a_1) \sqrt{N}.
\end{equation}
Taking $\varepsilon = \min\{ 1, t/a_1 \}$, we see that for each $x \in S^{n-1}$ satisfying
$\euclidnorm{\Gamma x} \leq t \sqrt{N}$ there is a corresponding $y \in \mcal{N}$ such that
$\euclidnorm{x - y} \leq \varepsilon$ and $\euclidnorm{\Gamma y} \leq 2t \sqrt{N}$.
Hence, using the union bound, setting $t = b_1/2$ and using Proposition \ref{propo-tall}, 
one has
%
%
$$
\PP (E \cap \Omega_{0})
\leq \sum_{y \in \mcal{N}} \PP \Bigl( \euclidnorm{\Gamma y} \leq 2t \sqrt{N} \Bigr)
\leq \cardinal{ \mcal{N} } e^{-b_{2}N}
\leq \Bigl( \frac{3}{\varepsilon} \Bigr)^n e^{-b_{2}N},
$$
where $b_{1}$ and $b_{2}$ are as in Proposition \ref{propo-tall}.
Thus
$$
 \PP (E \cap \Omega_{0}) \leq \exp \Bigl( -\frac{b_{2} N}{2} \Bigr)
$$
as long as
$$
 \Bigl( \frac{3}{\varepsilon} \Bigr)^n \leq \exp \Bigl( \frac{b_{2} N}{2} \Bigr).
$$
Bearing in mind that $N = (1 + \delta)n$, we can see that the last condition is satisfied if
\begin{equation}
 \delta \geq \delta_{0} := \max \left\{ \frac{2}{b_{2}} \ln
 \Bigl( \frac{6 a_1}{b_{1}} \Bigr), \; \frac{2}{b_{2}} \ln 3 \right\}.
\end{equation}
To finish, we use $\PP (E) \leq \PP (E \cap \Omega_{0}) + \PP (\Omega_{0}^c)$ 
with the estimate for $\PP (E \cap \Omega_{0})$ just obtained and the estimate 
$\PP (\Omega_{0}^c) \leq e^{-a_2 N}$ coming from condition (ii).
\hfill\qed

\section{``Almost square'' matrices (proof of Theorem~\ref{teo-almost})}

In this section we prove Theorem~\ref{teo-almost}. We will be using
all conditions (i) through (iv). The two key ingredients for the proof
of this theorem are Proposition \ref{propo-tall}
and Proposition \ref{propo-sbp}.
\medskip

\noindent
{\bf Proof of Theorem \ref{teo-almost}.} \quad
Passing to $r_0=\min\{3, r\}$ we may
assume without loss of generality that  $r\leq 3$.

Consider the event
$$
E := \{ \omega \st \exists x \in S^{n-1} 
\hspace{1mm} \text{ s.t. } \hspace{1mm} \euclidnorm{\Gamma x} \leq t \sqrt{N} \}.
$$
By equivalence \eqref{smallest-sv-equivalence} we are to estimate $\PP (E)$ 
with an appropriate value of $t$ (which will be specified later).

We split the set $E$ into two sets $E_{C}$ and $E_{I}$ defined as follows:
\begin{align*}
E_{C} &= \{ \omega \st \exists x \in S^{n-1} \cap Comp(m,\rho) 
\hspace{1mm} \text{ s.t. } \hspace{1mm} \euclidnorm{\Gamma x} \leq t \sqrt{N} \},\\
E_{I} &= \{ \omega \st \exists x \in S^{n-1} \cap Incomp(m,\rho) 
\hspace{1mm} \text{ s.t. } \hspace{1mm} \euclidnorm{\Gamma x} \leq t \sqrt{N} \},
\end{align*}
where $m \leq n$ and $\rho \in (0,1)$ will be specified later.

Define $\Omega_{0} := \{ \omega \st \norm{\Gamma} \leq a_1 \sqrt{N} \}$.
We will estimate $\PP (E)$ using the union bound in the inclusion
\begin{equation} \label{proof-teo-almost-1}
E \subset (E_{C} \cap \Omega_{0}) \cup (E_{I} \cap \Omega_{0}) \cup \Omega_{0}^c.
\end{equation}

Our proof will require that $t \leq 1$
(which will be satisfied once we choose $t$, see \eqref{proof-teo-almost-16} below);
and furthermore that $t$ and $\rho$ satisfy
\begin{equation} \label{proof-teo-almost-2}
\frac{2t}{a_1} \leq \rho \leq \frac{1}{4}.
\end{equation}

\noindent
\underline{\em Case 1: Probability of $E_{C} \cap \Omega_{0}$}.\quad
We work on the set $Comp(m,\rho)$, where $m \leq n$ and $\rho \in (0,1)$ 
will be specified later.

Given $x \in S^{n-1} \cap Comp(m,\rho)$, choose $y \in Sparse(m)$ so that 
$\euclidnorm{y - x} \leq \rho$.
It is clear that we may choose such a $y$ in $B_2^n$ 
(and thus  $1 - \rho \leq \euclidnorm{y} \leq 1$).
Note that on $\Omega_0$ we have $\norm{\Gamma} \leq a_1 \sqrt{N}$.
Thus if $x$ satisfies $\euclidnorm{\Gamma x} \leq t \sqrt{N}$ then
\begin{equation*}
\euclidnorm{\Gamma y}
\leq \euclidnorm{\Gamma x} + \norm{\Gamma}\cdot\euclidnorm{y - x}
\leq t \sqrt{N} + a_1 \rho \sqrt{N}  = (t + a_1 \rho) \sqrt{N}.
\end{equation*}

Let $\mcal{N}$ be a $\rho$-net in the set $B_2^n \cap Sparse(m)$. 
We may choose such a net with cardinality
\begin{displaymath}
\cardinal{\mcal{N}} \leq  \binom{n}{m} \Bigl( \frac{3}{\rho} \Bigr)^m
\leq \Bigl( \frac{en}{m} \Bigr)^m \Bigl( \frac{3}{\rho} \Bigr)^m
= \Bigl( \frac{3en}{\rho m} \Bigr)^m.
\end{displaymath}

For $y \in B_2^n \cap Sparse(m)$  chosen above, let $v \in \mcal{N}$ be such that 
$\euclidnorm{v - y} \leq \rho$. We observe that, by \eqref{proof-teo-almost-2},
$$
\euclidnorm{v} \geq \euclidnorm{y} - \rho \geq 1 - 2\rho \geq \frac{1}{2},
$$
and, by another use of \eqref{proof-teo-almost-2},
\begin{align*}
\euclidnorm{\Gamma v}
&\leq \euclidnorm{\Gamma y} + \norm{\Gamma}\cdot\euclidnorm{v - y} 
\leq (t + a_1 \rho) \sqrt{N} + \rho a_1 \sqrt{N}\\
&= (t + 2 a_1 \rho) \sqrt{N}
\leq \frac{5 a_1 \rho}{2} \sqrt{N}
\leq 5 a_1 \rho \sqrt{N} \euclidnorm{v}.
\end{align*}
Hence
\begin{equation}
\PP ( E_{C} \cap \Omega_{0} ) \leq
\PP \Bigl( \exists v \in \mcal{N} \text{ s.t. } 
   \euclidnorm{\Gamma v} \leq 5 a_1 \rho \sqrt{N} \euclidnorm{v} \Bigr) 
\leq \sum_{v \in \mcal{N}} \PP \Bigl( \euclidnorm{\Gamma v} \leq 5 a_1 \rho \sqrt{N} \euclidnorm{v} \Bigr). 
\label{proof-teo-almost-3}
\end{equation}
Using Proposition \ref{propo-tall}, we obtain
\begin{displaymath}
\PP \Bigl( \euclidnorm{\Gamma v} \leq 5 a_1 \rho \sqrt{N} \euclidnorm{v} \Bigr)
\leq e^{-b_2 N},
\end{displaymath}
provided that
\begin{equation} \label{proof-teo-almost-4}
5 a_1 \rho \leq b_{1}.
\end{equation}
We choose
\begin{equation} \label{proof-teo-almost-5}
\rho := \min\left\{ \frac{1}{4} \, , \, \frac{b_{1}}{5 a_1} \right\}
\end{equation}
so that both \eqref{proof-teo-almost-4} and the right hand side of \eqref{proof-teo-almost-2} 
are true. Now, from \eqref{proof-teo-almost-3}, we have
$$
\PP ( E_{C} \cap \Omega_{0} )
\leq \cardinal{ \mcal{N} } e^{-b_2 N}
\leq \Bigl( \frac{3en}{\rho m} \Bigr)^m e^{-b_2 N}.
$$
Thus, if
\begin{equation} \label{proof-teo-almost-6}
m \ln \Bigl( \frac{3en}{\rho m} \Bigr) \leq \frac{b_{2} N}{2}
\end{equation}
then
\begin{equation} \label{proof-teo-almost-7}
\PP ( E_{C} \cap \Omega_{0} ) \leq e^{-\frac{b_2 N}{2} }.
\end{equation}
Writing $m = \gamma n$, we see that inequality \eqref{proof-teo-almost-6} is satisfied if
\begin{displaymath}
\gamma \ln \Bigl( \frac{3e}{\rho \gamma} \Bigr) \leq \frac{b_{2}}{2},
\end{displaymath}
so we choose
\begin{equation} \label{proof-teo-almost-8}
\gamma = \frac{ b_2 }{ 4 \ln \bigl( \frac{6 e}{\rho b_2} \bigr) }.
\end{equation}

\noindent
\underline{\em Case 2: Probability of $E_{I} \cap \Omega_{0}$}.\quad
We work on the set $Incomp(m, \rho)$, where $\rho$ is defined in \eqref{proof-teo-almost-5}
and $m = \gamma n$ with $\gamma$ chosen in \eqref{proof-teo-almost-8}.

For convenience we set $a := t^{1/(r-2)}/a_1$. Since $t \leq 1$ and in view of
\eqref{proof-teo-almost-2}, we observe that $a \leq \rho/2$. Recall also that
that on $\Omega _0$ we have $\norm{\Gamma} \leq a_1 \sqrt{N}$.

Let $\mcal{N}$ be an $a$-net of $S^{n-1}$ with cardinality $\cardinal{\mcal{N}} \leq (3/a)^n$.
Let $x \in S^{n-1} \cap Incomp(m, \rho)$ be such that $\euclidnorm{ \Gamma x } \leq t \sqrt{N}$.
Recall that by (\ref{incompressible-equivalence}) one has
$\euclidnorm{ P_{\sigma} x} \geq \frac{\rho}{2}$ for every
$\sigma \subset \{ 1, \ldots, n \}$  with $\cardinal{\sigma^c}\leq m$.
Then there is $v \in \mcal{N}$ such that $\euclidnorm{ \Gamma v } \leq 2t \sqrt{N}$
and with the additional property $\euclidnorm{ P_{\sigma} v } \geq \frac{\rho}{2}$
for each $\sigma \subset \{ 1, \ldots, n \}$  with $\cardinal{\sigma^c} \leq m$.
Indeed, choosing $v \in \mcal{N}$ such that $\euclidnorm{x - v} \leq a$ and using
$a_1 a = t^{1/(r-2)} \leq t$ (which holds by the choice of $a$), we have
\begin{equation*}
\euclidnorm{\Gamma v}
\leq \euclidnorm{\Gamma x} + \norm{\Gamma} \cdot \euclidnorm{v - x}
\leq t \sqrt{N} + a_1 \sqrt{N} a  \leq 2t \sqrt{N}
\end{equation*}
and
\begin{equation*}
\euclidnorm{P_{\sigma} v}
\geq \euclidnorm{P_{\sigma} x} - \euclidnorm{P_{\sigma} (v - x)}
\geq \rho - a \geq \frac{\rho}{2},
\end{equation*}
where we used the condition $2a \leq 2t/a_1 \leq \rho$, required in \eqref{proof-teo-almost-2}.

Denote by $\mcal{A}$ the set of all $v \in \mcal{N}$ with the property that for each set
$\sigma \subset \{ 1, \ldots, n \}$  with $\cardinal{\sigma^c} \leq m$ we have 
$\euclidnorm{ P_{\sigma} v } \geq \frac{\rho}{2}$.
Then
\begin{equation} \label{proof-teo-almost-9}
\PP ( E_{I} \cap \Omega_{0} ) \leq
\PP \Bigl( \exists v \in \mcal{A} \st \euclidnorm{ \Gamma v } \leq 2t \sqrt{N} \Bigr).
\end{equation}
Now, for each fixed $v = (v_i) \in \mcal{A}$ we have
\begin{align}
\PP \Bigl( \euclidnorm{ \Gamma v }^2 \leq 4t^2N \Bigr)
&= \PP \Bigl( N - \frac{ 1 }{ 4t^2 } \euclidnorm{ \Gamma v }^2 \geq 0 \Bigr) \notag\\
&\leq \EE \exp \Bigl\{ N - \frac{ 1 }{ 4t^2 } \euclidnorm{ \Gamma v }^2 \Bigr\} \notag\\
&= e^N \EE \exp \Bigl\{ - \frac{ 1 }{ 4t^2 } \sum_{j=1}^{N} \Bigabs{ \sum_{i=1}^{n} \xi_{ji} v_i }^2 \Bigr\} \notag\\
&= e^N \prod_{j=1}^{N} \EE \exp \Bigl\{ - \frac{ 1 }{ 4t^2 } \Bigabs{ \sum_{i=1}^{n} \xi_{ji} v_i }^2 \Bigr\}, 
\label{proof-teo-almost-10}
\end{align}
and our goal is to make this last expression small.
To estimate the expectations we use the distribution formula:
\begin{align}
\EE \exp \Bigl\{ - \frac{ 1 }{ 4t^2 } \Bigabs{ \sum_{i=1}^{n} \xi_{ji} v_i }^2 \Bigr\}
&= \int_{0}^{1}
\PP \biggl( \exp \Bigl\{ - \frac{ 1 }{ 4t^2 } \Bigabs{ \sum_{i=1}^{n} \xi_{ji} v_i }^2 \Bigr\} > s \biggr) ds \notag\\
&= \int_{0}^{\infty} u e^{-u^2/2} \PP \biggl( \Bigabs{ \sum_{i=1}^{n} \xi_{ji}v_i } < \sqrt{2} tu \biggr) du. 
\label{proof-teo-almost-11}
\end{align}
It is now apparent that we need to estimate the quantities
\begin{displaymath}
f_j(\lambda) := \PP \biggl( \Bigabs{ \sum_{i=1}^{n} \xi_{ji} v_i } < \lambda \biggr),
\quad j\leq N.
\end{displaymath}
To this end, note that for each row $j \in \{ 1, \ldots, N \}$ there exists $\sigma_j \subset \{ 1, \ldots, n \}$
with cardinality $\cardinal{\sigma_j} \geq a_4 n$ such that $\EE \xi_{ji}^2 \geq 1$ for all $i \in \sigma_j$
(this is condition (iv)). Also, for each fixed $v$, set
\begin{displaymath}
\sigma_{v} := \{ i \st \abs{v_i} > a \}.
\end{displaymath}
Since $v \in S^{n-1}$ we have $\cardinal{ \sigma_{v} } \leq 1/a^2$.

Set $\overline{\sigma}_j = \sigma_j \setminus \sigma_v$, and note that
\begin{displaymath}
\cardinal{ \overline{\sigma}_j } \geq a_4 n - \frac{1}{a^2}.
\end{displaymath}
It follows that $\cardinal{ \overline{\sigma}_j^c } \leq (1 - a_4) n + \frac{1}{a^2}$,
so to have $\cardinal{ \overline{\sigma}_j^c } \leq m$ it suffices to require
\begin{equation} \label{proof-teo-almost-12}
(1 - a_4) n + \frac{1}{a^2} \leq m.
\end{equation}
Note that (\ref{proof-teo-almost-12}), in particular, implies $1/a^2 \leq a_4 n \leq n$.
Recall that $m = \gamma n$, where $\gamma$ was chosen in \eqref{proof-teo-almost-8}.
Then inequality (\ref{proof-teo-almost-12}) is satisfied if $a_4 > 1 - \gamma$
(which is the condition on $\gamma$ in our Theorem) and
\begin{equation} \label{proof-teo-almost-13}
t \geq \biggl( \frac{a_1}{\sqrt{(\gamma + a_4 - 1) n}} \biggr)^{r-2}.
\end{equation}
Now, since $\cardinal{ \overline{\sigma}_j^c } \leq m$, we have 
$\euclidnorm{ P_{\overline{\sigma}_j} v } \geq \rho/2 $,
and hence
\begin{displaymath}
A_j^2 := \EE \sum_{i \in \overline{\sigma}_j} \xi_{ji}^2 v_i^2 \geq \frac{\rho^2}{4}
\end{displaymath}
(where we have used the property $\EE \xi_{ji}^2 \geq 1$ for $i \in \sigma_j$).
Consequently, using Proposition \ref{propo-sbp}, and keeping in mind $\abs{v_i} \leq a$ 
for $i \in \overline{\sigma}_j$,
we get
$$
f_j( \lambda )
\leq c \Bigl( \frac{\lambda}{\rho} + \frac{ \mu^r }{\rho^r} \norm{ P_{\overline{\sigma}_j} v }_r^r \Bigr) 
\leq c \Bigl( \frac{\lambda}{\rho} + \frac{ \mu^r }{\rho^r} \norm{ P_{\overline{\sigma}_j} v }_{\infty}^{r-2} \cdot
\euclidnorm{ P_{\overline{\sigma}_j} v }^2 \Bigr) 
\leq c \Bigl( \frac{\lambda}{\rho} + \frac{ \mu^r a^{r-2} }{\rho^r} \Bigr)
$$
for some absolute constant $c \geq 1$.
Then, continuing from \eqref{proof-teo-almost-11} we have
\begin{align*}
\EE \exp \Bigl\{ - \frac{ 1 }{ 4t^2 } \Bigabs{ \sum_{i=1}^{n} \xi_{ji} v_i }^2 \Bigr\}
&= \int_{0}^{\infty} u e^{-u^2/2} f_j( \sqrt{2} tu ) du \\
&\leq c \int_{0}^{\infty} u e^{-u^2/2} \Bigl( \frac{\sqrt{2} tu}{\rho} + \frac{ \mu^r a^{r-2} }{\rho^r} \Bigr) du \\
&= \frac{c \sqrt{2} t }{\rho} \int_{0}^{\infty} u^2 e^{-u^2/2} du +
     \frac{c \mu^r a^{r-2} }{\rho^r} \int_{0}^{\infty} u e^{-u^2/2} du \\
&= \frac{c \sqrt{\pi} t }{\rho} + \frac{c \mu^r t}{\rho^r a_1^{r-2}} = c_3 t,
\end{align*}
where
\begin{equation} \label{proof-teo-almost-14}
c_3 := c \Bigl( \frac{\sqrt{\pi}}{\rho} + \frac{\mu^r}{\rho^r a_1^{r-2}} \Bigr).
\end{equation}
Therefore, from \eqref{proof-teo-almost-10}, we get (for each fixed $v \in \mcal{A}$)
\begin{displaymath}
\PP \Bigl( \euclidnorm{ \Gamma v } \leq 2t \sqrt{N} \Bigr) \leq e^N ( c_3 t )^N = (c_3 e t)^N,
\end{displaymath}
and from this, in \eqref{proof-teo-almost-9} we get
$$
\PP ( E_{I} \cap \Omega_{0} )
\leq \cardinal{ \mcal{A} } (c_3 e t)^N
\leq \Bigl( \frac{3}{a} \Bigr)^n (c_3 e t)^N
= \Bigl( \frac{3 a_1}{t} \Bigr)^n (c_3 e t)^N.
$$
Then we can make
\begin{equation} \label{proof-teo-almost-15}
\PP ( E_{I} \cap \Omega_{0} ) \leq e^{-N}
\end{equation}
provided that
\begin{equation} \label{proof-teo-almost-16}
t \leq  \frac{1}{c_3 e^2} \  \Bigl( \frac{1}{3 a_1 c_3 e^2} \Bigr)^{1/\delta}.
\end{equation}

Choose $t$ to satisfy equality in \eqref{proof-teo-almost-16}. Note
\begin{equation}
  \frac{t}{a_1} \leq \frac{1}{c_3 e^2 a_1} \leq \frac{\rho}{c e^2 \sqrt{\pi}} ,
\end{equation}
so the left hand side of \eqref{proof-teo-almost-2} holds.
Finally note that \eqref{proof-teo-almost-13} is satisfied whenever
$$
\delta \geq 
\frac{ \frac{2}{r-2} \ln (3 a_1 c_3 e^2) }{ \ln \Bigl( \frac{ (\gamma + a_4 - 1) n }{ a_1^2 (c_3 e^2)^{2/(r-2)} } \Bigr) }
\; =: \; \frac{\tilde{c}_1}{ \ln (\tilde{c}_2 n) }.
$$
To finish, we take probabilities in \eqref{proof-teo-almost-1} and we use the estimates for
$\PP ( E_{C} \cap \Omega_{0} )$ and $\PP ( E_{I} \cap \Omega_{0} )$
we have found in \eqref{proof-teo-almost-7} and \eqref{proof-teo-almost-15}, respectively, 
combined with the estimate $\PP ( \Omega_0^c ) \leq e^{-a_2 N}$ coming from condition (ii).
This shows that, with the chosen $t$, we have $\PP (E) \leq e^{-b_2 N/2} + e^{-N} + e^{-a_2 N}$,
which completes the proof.
\qed

\section{Square matrices (proof of Theorem~\ref{teo-square})}

In this section our goal is to prove Theorem~\ref{teo-square}. We are going to use
two lemmas from \cite{RV1}. The first one is \cite[Lemma 3.5]{RV1}.
Note that the proof given there works for any random matrix.

\begin{lemma}
\label{lema-invert-via-distance}
Let  $\Gamma$ be any random matrix of size $m \times n$.
Let $X_1, \ldots, X_n$ denote the columns of $\Gamma$ and let $H_k$ denote the span of
all column vectors except the $k$-th. Then for every $\gamma, \rho \in (0,1)$ and every
$\varepsilon > 0$ one has
\begin{displaymath}
 \PP \Bigl( \inf_{x \in F} \euclidnorm{\Gamma x} \leq \varepsilon \rho n^{-1/2} \Bigr)
 \leq \frac{1}{\gamma n} \sum_{k = 1}^{n} \PP \bigl( \dist(X_k,H_k) < \varepsilon \bigr),
\end{displaymath}
where $F= S^{n-1} \cap Incomp(\gamma n,\rho)$.
\end{lemma}

The next lemma is similar to \cite[Lemma 3.8]{RV1}.
To prove it one would repeat the proof of that lemma,
replacing \cite[Lemma 3.7]{RV1} used there with our Lemma~\ref{sbp-incompressible}.

\begin{lemma}
\label{lema-weak-distance-bound}
Let $r\in (2, 3]$ and $\Gamma$ be a random matrix as in Theorem \ref{teo-square}.
Let $X_1, \ldots, X_n$ denote its column vectors, and consider the subspace
$H_n = \linspan (X_1, \ldots, X_{n-1})$. 
Then there exists a positive constant $\gamma_0 =\gamma_0(r,\mu, a_1, a_3)< 1$ 
such that if $a_4> 1 - \gamma_0$ then 
 for every $\varepsilon \geq 0$ one has
\begin{displaymath}
  \PP \Bigl( \dist(X_n,H_n) < \varepsilon \quad \mbox{ and } \quad \norm{\Gamma} \leq
  a_1 n^{1/2} \Bigr) \leq c (\varepsilon n^{\frac{3-r}{2}} + \mu^r n^{\frac{2-r}{2}}),
\end{displaymath}
where $c$ depends on $r$, $\mu$, $a_1$, $a_3$, and $a_4$.
\end{lemma}

Now we are ready for the proof of Theorem \ref{teo-square}.

\medskip

\noindent
{\bf Proof of Theorem \ref{teo-square}.} \quad
Without loss of generality we assume $\varepsilon \leq a_1/2$
(otherwise choose $C = 2/a_1$ and we are done). We also assume that
$r\leq 3$ (otherwise we pass to $r_0=\min\{3, r\}$).

Consider the event
$$
E := \{ \omega \st \exists x \in S^{n-1} \hspace{1mm} \text{ s.t. } \hspace{1mm}
\euclidnorm{\Gamma x} \leq t n^{-1/2} \}.
$$
By equivalence \eqref{smallest-sv-equivalence} we are to estimate $\PP (E)$
with an appropriate value of $t$ (which will be specified later).

As in the proof of Theorem \ref{teo-almost},
we split the set $E$ into the sets $E_{C}$ and $E_{I}$ defined as follows:
\begin{align*}
E_{C} &= \{ \omega \st \exists x \in S^{n-1} \cap Comp(m, \rho) 
\hspace{1mm} \text{ s.t. } \hspace{1mm} \euclidnorm{\Gamma x} \leq t n^{-1/2} \},\\
E_{I} &= \{ \omega \st \exists x \in S^{n-1} \cap Incomp(m, \rho) 
\hspace{1mm} \text{ s.t. } \hspace{1mm} \euclidnorm{\Gamma x} \leq t n^{-1/2} \}.
\end{align*}

Define $\Omega_{0} := \{ \omega \st \norm{\Gamma} \leq a_1 \sqrt{n} \}$.
We will estimate $\PP (E)$ using the union bound in the inclusion
\begin{equation} \label{proof-teo-square-1}
E \subset (E_{C} \cap \Omega_{0}) \cup E_{I} \cup \Omega_{0}^c.
\end{equation}

\noindent
\underline{\em Case 1: Probability of $E_{C} \cap \Omega_{0}$}.\quad
The proof of this case is almost line to line repetition of the corresponding
proof in Theorem~\ref{teo-almost} (see Case~1 there).
Let $m \leq n$ and $\rho \in (0,1)$ be specified later.
Using approximation argument and the union bound as
in the proof of  Case 1 in Theorem~\ref{teo-almost},
and choosing
\begin{equation} \label{proof-teo-square-2}
\rho := \min\left\{ \frac{1}{4} \, , \, \frac{b_1}{5 a_1} \right\}, \hspace{1cm}
\gamma := \frac{b_2}{ 4 \ln \bigl( \frac{6 e}{\rho b_2} \bigr) }, \hspace{1cm}
m = \gamma n,
\end{equation}
we obtain
\begin{equation} \label{proof-teo-square-3}
\PP ( E_{C} \cap \Omega_{0} ) \leq e^{-b_2 n/2},
\end{equation}
provided that
\begin{equation} \label{proof-teo-square-4}
\frac{2t}{a_1} \leq \rho.
\end{equation}

\noindent
\underline{\em Case 2: Probability of $E_{I}$}.\quad
We work on the set $Incomp(m, \rho)$, where $m = \gamma n$ and $\gamma$, $\rho$ 
chosen in \eqref{proof-teo-square-2}.

Using Lemma \ref{lema-invert-via-distance} with $\varepsilon = t/\rho$, 
and also applying Lemma \ref{lema-weak-distance-bound},
we get 
\begin{align}\PP ( E_{I} )
& \leq \frac{1}{\gamma n} \sum_{k = 1}^{n} \PP \bigl( \dist(X_k,H_k) < t/\rho \bigr) \notag \\
& \leq \frac{1}{\gamma n}
       \sum_{k = 1}^{n} \Bigl\{ \PP \bigl( \dist(X_k,H_k) < t/\rho \quad \& \quad \norm{\Gamma} \leq a_1 \sqrt{n} \bigr)
             + \PP \bigl(\norm{\Gamma} > a_1 \sqrt{n} \bigr) \Bigr\} \notag \\
& \leq \frac{1}{\gamma n} \sum_{k = 1}^{n} \Bigl\{ c (\varepsilon n^{\frac{3-r}{2}} + n^{\frac{2-r}{2}}) + e^{- a_2 n} \Bigr\} \notag \\
& \leq \frac{c}{\gamma} (\varepsilon n^{\frac{3-r}{2}} + n^{\frac{2-r}{2}}) + \frac{1}{\gamma} e^{- a_2 n}.
\label{proof-teo-square-5}
\end{align}

Also notice that our choice $t = \varepsilon \rho$ and our assumption $\varepsilon \leq a_1/2$
guarantee that $t$ satisfies \eqref{proof-teo-square-4}.

To finish the proof, we take probabilities in \eqref{proof-teo-square-1}, and we use 
the estimates for $\PP ( E_{C} \cap \Omega_{0} )$ and for $\PP ( E_{I} )$ obtained in 
\eqref{proof-teo-square-3} and \eqref{proof-teo-square-5},
respectively, combined with the estimate $\PP ( \Omega_0^c ) \leq e^{-a_2 n}$ 
coming from condition (ii). This way we obtain
\begin{align*}
 \PP(E) \leq e^{-b_2 n/2} + \frac{c}{\gamma} (\varepsilon n^{\frac{3-r}{2}} + n^{\frac{2-r}{2}})
  + \frac{1}{\gamma} e^{- a_2 n} + e^{-a_2 n} \leq C (\varepsilon n^{\frac{3-r}{2}} + n^{\frac{2-r}{2}})
\end{align*}
for a suitable constant $C$.
\qed

\section{Proof of Proposition \ref{propo-tall}}

Take an arbitrary $x = (x_1, \ldots, x_n) \in \R^n$ with $\euclidnorm{x} = 1$.
For $a > 0$ (a parameter whose value will be specified later),
define a set of ``good'' rows as follows:
$$
J = J(a) = \biggl\{ j\in \{ 1, \ldots, N\} \st \EE \sum_{i=1}^{n} \xi_{ji}^2 x_i^2 \geq a \biggr\}.
$$
Suppose that the cardinality of set $J$ is $\cardinal{ J } = \alpha N$
for some $\alpha \in [0,1]$.
Note that for each index $j = 1, \ldots, N$ we have
$$
\EE \sum_{i=1}^{n} \xi_{ji}^2 x_i^2
\leq \max_{1\leq i\leq n} \EE \xi_{ji}^2
\leq \max_{1\leq i\leq n} (\EE \xi_{ji}^r)^{2/r}
\leq \mu^2.
$$
Then on one hand we have
\begin{align*}
\sum_{j=1}^{N} \biggl( \EE \sum_{i=1}^{n} \xi_{ji}^2 x_i^2 \biggr)
&= \sum_{j \in J} \biggl( \EE \sum_{i=1}^{n} \xi_{ji}^2 x_i^2 \biggr) +
\sum_{j \in J^c} \biggl( \EE \sum_{i=1}^{n} \xi_{ji}^2 x_i^2 \biggr) \\[.2cm]
&\leq \mu^2 \alpha N + a(1-\alpha)N,
\end{align*}
while on the other hand, using condition (iii),
$$
\sum_{j=1}^{N} \biggl( \EE \sum_{i=1}^{n} \xi_{ji}^2 x_i^2 \biggr)
= \sum_{i=1}^{n} \biggl( \EE \sum_{j=1}^{N} \xi_{ji}^2 \biggr) x_i^2
\geq \sum_{i=1}^{n} a_3^2 N x_i^2 = a_3^2 N.
$$
Hence we have $\mu^2 \alpha N + a(1-\alpha)N \geq a_3^2 N$, so $\alpha$ satisfies
\begin{equation}
\alpha \geq \frac{a_3^2 - a}{\mu^2 - a}. \label{proof-propo-tall-1}
\end{equation}
Note that for each $j = 1,\ldots,N$, the $j$-th entry of $\Gamma x$ is
$(\Gamma x)_j = \sum_{i=1}^{n} \xi_{ji} x_i$.
Define $f_j := \bigabs{\sum_{i=1}^{n} \xi_{ji} x_i}$, so
$$
\euclidnorm{\Gamma x}^2 = \sum_{j=1}^{N} f_j^2.
$$
Clearly $f_1,\ldots,f_N$ are independent.
For any $t, \tau > 0$ we have
\begin{align}
\PP \bigl( \euclidnorm{\Gamma x}^2 \leq t^2 N \bigr)
&= \PP \biggl( \sum_{j=1}^{N} f_j^2 \leq t^2 N \biggr) 
= \PP \biggl( \tau N - \frac{\tau}{t^2} \sum_{j=1}^{N} f_j^2 \geq 0 \biggr) \notag\\
&\leq \EE \exp \biggl( \tau N - \frac{\tau}{t^2} \sum_{j=1}^{N} f_j^2 \biggr) 
= e^{\tau N} \prod_{j=1}^{N} \EE \exp \biggl( -\frac{\tau f_j^2}{t^2} \biggr). 
\label{proof-propo-tall-2}
\end{align}
%
%
From Lemma \ref{lema-deviation} we know that for every $j=1,\ldots,N$,
\begin{equation} \label{probmy}
\PP (f_j > \lambda) \geq
\biggl( \frac{ [ \EE \sum_{i=1}^{n} \xi_{ji}^2 x_i^2 - \lambda^2]_+ }{ 8\mu^2 } \biggr)^{r/(r - 2)} =: \beta_j(r),
\end{equation}
Note that for every $j \in J$ one has
\begin{equation}
\beta_j \geq \left( \frac{ [a-\lambda^2]_+ }{ 8\mu^2 } \right)^{r/(r - 2)}.
\label{proof-propo-tall-3}
\end{equation}
For arbitrary $t > 0$, $\eta > 0$ and $\lambda > 0$, set $\tau := \frac{\eta t^2}{\lambda ^2}$.
For each $j = 1, \ldots, N$ we have
\begin{align*}
\EE \exp \Bigl( -\frac{\tau f_j^2}{t^2} \Bigr)
&= \int_{0}^{1} \PP \biggl( \exp \Bigl( -\frac{\eta f_j^2}{\lambda^2} \Bigr) > s \biggr) ds\\
&= \int_{0}^{e^{-\eta}} \!\!\! \PP \biggl( \exp \Bigl( \frac{\eta f_j^2}{\lambda^2} \Bigr) < \frac{1}{s} \biggr) ds +
\int_{e^{-\eta}}^{1} \!\!\! \PP \biggl( \exp \Bigl( \frac{\eta f_j^2}{\lambda^2} \Bigr) < \frac{1}{s} \biggr) ds\\
&\leq e^{-\eta} + \PP (f_j < \lambda) ( 1 - e^{-\eta} ) .
\end{align*}
Choosing $\eta = \ln 2$ and applying (\ref{probmy}), we obtain
$$
 \EE \exp \Bigl( -\frac{\tau f_j^2}{t^2} \Bigr) \leq
 e^{-\eta} + (1 - \beta_j(r))( 1 - e^{-\eta} ) = 1 - \frac{\beta_j(r)}{2}
 \leq \exp\left( - \frac{\beta_j(r)}{2} \right) .
$$
Since $\tau < \frac{t^2}{\lambda ^2}$, inequality \eqref{proof-propo-tall-2} implies
\begin{equation}\label{proof-propo-tall-4}
  \PP \bigl( \euclidnorm{\Gamma x}^2 \leq t^2 N \bigr)
  \leq e^{\tau N} \prod_{j=1}^{N} e^{-\beta_j(r)/2}
 \leq e^{(t^2/\lambda^2)N} \prod_{j \in J} e^{-\beta_j(r)/2}.
\end{equation}
Taking $a = a_3^2/2$ and  $\lambda = a_3/2$ and using \eqref{proof-propo-tall-3}
we observe that for every $j \in J$
we have $\beta_j \geq \bigl( \frac{a_3^2}{32\mu^2} \bigr)^{r/(r-2)}$.
Also note this choice of $a$ and \eqref{proof-propo-tall-1} imply $\alpha \geq a_3^2/(2\mu^2)$.
Now let
$$
 t^2 := \frac{a_3^{4}}{2^{5}\mu^2} \Bigl( \frac{a_3^{2}}{2^{5}\mu^2} \Bigr)^{r/(r - 2)}.
$$
Then continuing from \eqref{proof-propo-tall-4}
we obtain
$$
\PP \biggl( \euclidnorm{\Gamma x}^2 \leq \frac{a_3^{4}}{2^{5}\mu^2} \Bigl( \frac{a_3^{2}}{2^{5}\mu^2} \Bigr)^{r/(r - 2)} N \biggr)
\leq \exp \biggl\{ - \frac{a_3^{2}}{2^{3}\mu^2} \Bigl( \frac{a_3^{2}}{2^{5}\mu^2} \Bigr)^{r/(r - 2)} N \biggr\}.
$$
This completes the proof.
\qed

\bigskip


\smallskip

\noindent
A. E. Litvak,
{\small Dept. of Math. and Stat. Sciences},
{\small University of Alberta}, {\small Edmonton, Alberta T6G 2G1, Canada},
{\small \tt  aelitvak@gmail.com}

\bigskip

\noindent
O.  Rivasplata,
{\small Dept. of Math. and Stat. Sciences},
{\small University of Alberta}, {\small Edmonton, Alberta T6G 2G1, Canada},
{\small \tt  orivasplata@ualberta.ca}

\end{document}